\documentclass{amsart}
\usepackage{amssymb,amsmath}
\newcommand{\RR}{{\mathbb R}}

\newcommand{\e}{\varepsilon}

\newcommand{\Lap}{\Delta}
\newcommand{\de}{\delta}
\newcommand{\del}{\partial}
\newcommand{\om}{\omega}

\newcommand{\si}{\sigma}

\newcommand{{\loc}}{{\ell\mathrm oc}}

\def\meanint{{\diagup\hskip -.42cm\int}}

\newtheorem{theorem}{Theorem}
\newtheorem{lemma}{Lemma}
\newtheorem{proposition}{Proposition}
\newtheorem{corollary}{Corollary}

\begin{document}

\title[Fundamental Solution with Singular Drift]{The Fundamental Solution of an Elliptic Equation with Singular Drift}
\author{Vladimir Maz'ya}
\address{Department of Mathematics, Link\"oping University, SE-581 83 Link\"oping, Sweden 
and RUDN, 6 Miklukho-Maklay St, Moscow, 117198, Russia
}
\email{vladimir.mazya@liu.se}

\author{ Robert McOwen}
\address{Department of Mathematics, Northeastern University, Boston, MA 02115}
\email{r.mcowen@northeastern.edu}
\date{August 31, 2022}

\begin{abstract}
For $n\geq 3$, we study the existence and asymptotic properties of the fundamental solution for elliptic operators in nondivergence form,
${\mathcal L}(x,\partial_x)=a_{ij}(x)\partial_i\partial_j+b_k(x)\partial_k$, where the $a_{ij}$ have modulus of continuity $\om(r)$ satisfying the square-Dini condition and the $b_k$ are allowed mild singularities of order $r^{-1}\om(r)$. A singular integral is introduced that controls the existence of the fundamental solution. We give examples that show the singular drift $b_k\partial_k$ may  act as a perturbation that does not dramatically change the fundamental solution of ${\mathcal L}^o=a_{ij}\partial_i\partial_j$, or it could change an operator ${\mathcal L}^o$ that does not have a fundamental solution to one that does. 
\end{abstract}


\maketitle

\addtocounter{section}{-1}
\section{Introduction}\label{sec:0}

Consider a 2nd-order operator in nondivergence form with a 1st-order term
\begin{equation}\label{L-nondivergence}
{\mathcal L}(x,\partial_x)u= a_{ij}(x)\,\partial_i \partial_j u+b_k(x)\,\partial_k u, \quad\hbox{for $x\in U$},
 \end{equation}
 where the coefficients $a_{ij}=a_{ji}$ and $b_k$ are measurable, real-valued functions in an open set $U\subset\RR^n$, $n\geq 3$; here and throughout the paper we use the summation convention for repeated indices. 
 We assume ${\mathcal L}(x,\partial_x)$ is pointwise elliptic, i.e.\ the matrix $a_{ij}(x)$ is positive definite for each $x\in U$. 
 Our objective is to study the existence and asymptotics of its  {\it fundamental solution}, i.e.\ a function $F(x,y)$ satisfying $F(x,\cdot)\in L^1_{\loc}(U)$ for each $x\in U$ and 
 \begin{equation}\label{LF=delta}
-{\mathcal L}(x,\partial_x) F(x,y)=\delta(x-y)\quad\hbox{for}\ x,y\in U,
 \end{equation}
 in a distributional sense that needs to be made clear. In the classical case that the coefficients are sufficiently smooth (e.g.\ H\"older continuous), the existence of a fundamental solution and its asymptotic properties are well-known: cf.\  Miranda \cite{M}.
  
 In \cite{MM1} we studied this problem for ${\mathcal L}^o(x,\partial_x)=a_{ij}(x)\,\partial_i\partial_j$ when the modulus of continuity $\om(r)$ for the $a_{ij}$ satisfies the {\it square-Dini condition}
  \begin{equation}\label{squareDini}
  \int_0^1 \om^2(r)\frac{dr}{r}<\infty.
   \end{equation}
   Notice that this assumption is weaker than $\lambda$-H\"older continuity since \eqref{squareDini} is satisfied by $\om(r)=r^\lambda$ for $0<\lambda<1$; \eqref{squareDini} is also weaker than the Dini-condition, $\int_0^1 r^{-1}\om(r)\,dr<\infty$. The hypotheses and conclusions of  \cite{MM1} are most easily stated if we fix $y=0$ and let $U=B_\e=\{x\in\RR^n:|x|<\e\}$ for $\e$ sufficiently small.  By an affine change of coordinates, we may arrange  $a_{ij}(0)=\delta_{ij}$, so we may assume
   \begin{equation}\label{aij-delta}
 \sup_{|x|=r} |a_{ij}(x)-\delta_{ij}|\leq \om(r)\quad\hbox{for}\ 0<r<\e.
 \end{equation}
It was found in \cite{MM1} that the existence of a solution of 
    \begin{equation}\label{LF=delta-y=0}
-{\mathcal L}^o(x,\partial_x) F(x)=\delta(x)\quad\hbox{in} \ B_\e,
 \end{equation}
for $\e$ sufficiently small depends on the behavior of the integral
\begin{equation}\label{def:I0}
I^o(r)=\frac{1}{|S^{n-1}|}\int_{r<|z|<\e} \left(\hbox{tr}({\bf A}_z)-n\frac{\left\langle {\bf A}_z z,z\right\rangle}{|z|^2}\right)\frac{dz}{|z|^n}
\end{equation}
as $r\to 0$; here ${\bf A}_x$ denotes the matrix $(a_{ij})$ and $\langle {\bf A}_z z,z\rangle=a_{ij}(z)z_iz_j$. If $I^o(0)=\lim_{r\to 0} I^o(r)$ exists and is finite, then we showed there is a solution of \eqref{LF=delta-y=0}.  However, if $I^o(r)\to -\infty$ as $r\to 0$, then we found a solution $F(x)$ of
    \begin{equation}\label{LF=0}
{\mathcal L}^o(x,\partial_x) F(x)=0\quad\hbox{in} \ B_\e,
 \end{equation}
 which has a singularity at $x=0$, but $F(x)=o(|x|^{2-n})$ as $|x|\to 0$ and $F(x)$ is not a solution of \eqref{LF=delta-y=0}; this violates the ``extended maximum principle'' of \cite{GS}.

In the present paper, we generalize the results of  \cite{MM1} by including the 1st-order term $b_k\partial_k$, which is called a {\it drift} term: cf.\ \cite{CZ}, \cite{ABMMZ}. 
If the $b_k$ are bounded in $U$, they will  not effect the existence of the fundamental solution, so we will allow the drift to be singular at $x=0$, but satisfy the condition:
\begin{equation}\label{est:b_k}
\sup_{|x|=r} |x|\,|b_k(x)|\leq c\,\om(r)\quad\hbox{for}\ 0<r<\e.
\end{equation}
(Here and throughout this paper, $c$ denotes a generic constant.)
We again construct an unbounded solution $Z(x)$ of
 \begin{equation}\label{LZ=0}
 {\mathcal L}(x,\partial_x)Z(x)=0 \quad\hbox{for}\ x\in B_\e \backslash \{0\},
  \end{equation}
 and then check to see whether $Z$ has the appropriate singular behavior as $|x|\to 0$ so that, in a distributional sense,
  \begin{equation}\label{LZ=C0*delta}
- {\mathcal L}(x,\partial_x)Z(x)=C_0\,\delta(x) \quad\hbox{for some constant $C_0$.}
  \end{equation}
  If so, then setting $F(x):=C_0^{-1}\,Z(x)$ defines a solution of
  \begin{equation}\label{-LF=delta}
  - {\mathcal L}(x,\partial_x)F(x)=\delta(x) \quad\hbox{for} \ x\in B_\e.
  \end{equation}
To determine whether this can be done, the quantity \eqref{def:I0} is generalized to 
\begin{equation}\label{def:I}
I(r)=\frac{1}{|S^{n-1}|}\int_{r<|z|<\e} \left(\hbox{tr}({\bf A}_z)-n\frac{\left\langle {\bf A}_z z,z\right\rangle}{|z|^2}+\langle{\bf B}_z,z\rangle\right)\frac{dz}{|z|^n},
\end{equation}
where ${\bf B}_z$ denotes the vector with components $b_k(z)$ and $\langle{\bf B}_z,z\rangle=b_k(z)\,z_k$. In general, $I(r)$ may have a singularity at $r=0$, but it is weaker than logarithmic: for any $\lambda>0$ there exists $c_\lambda>0$ such that
    \begin{equation}
    |I(r)|\leq \lambda \,|\log r|+c_\lambda\quad\hbox{for}\ 0<r<\e.
     \end{equation}
  
Before stating the main results of this paper, we need some notation and an additional assumption.
For an open set $U$, $1< p<\infty$, and integer $m=0,1,2$, let $W^{m,p}_\loc(U)$ denote the Sobolev space of functions whose derivatives up to order $m$ are locally $L^p$-integrable in $U$. The $L^p$-mean of $w\in L^{p}_\loc(\RR^n\backslash\{0\})$ on the annulus $A_r:=\{x: r<|x|<2r\}$ is defined by
\begin{subequations}
\begin{equation}
M_p(w,r):=\left(\meanint_{A_r} |w(x)|^p\,dx\right)^{1/p}.
\end{equation}
Similarly, we define
\begin{equation}
M_{1,p}(w,r)=rM_p(Dw,r)+M_p(w,r) \quad\hbox{for}\ w\in W^{1,p}_\loc(\RR^n\backslash\{0\}),
\end{equation}
\begin{equation}
M_{2,p}(w,r)=r^2 M_p(D^2w,r)+M_{1,p}(w,r) \quad\hbox{for}\ w\in W^{2,p}_\loc(\RR^n\backslash\{0\}).
\end{equation}
\end{subequations}
For $p=\infty$ we can analogously define $M_\infty(w,r)$,  $M_{1,\infty}(w,r)$ and $M_{2,\infty}(w,r)$.
The $L^p$-mean over annuli centered at $y\not=0$ will be denoted by $M_p(w,r;y)$ and similarly for $M_{1,p}(w,r;y)$ and $M_{2,p}(w,r;y)$.
We also define
   \begin{equation}\label{sigma(r)=}
   \sigma(r)=\int_0^r  \om^2(\rho)\frac{d\rho}{\rho},
   \end{equation}
   which satisfies $\sigma(r)\to 0$ as $r\to 0$ because of \eqref{squareDini}. Finally, as a modulus of continuity we want $\om(r)$ to be
   nondecreasing for $0<r<1$, but  we also assume for some $\kappa\in (0,1)$ that 
   \begin{equation}\label{om-monotone}
   \om(r)\,r^{-1+\kappa}\quad\hbox{is nonincreasing on $0<r<1$.}
   \end{equation}
   This is natural since we are interested in moduli of continuity that vanish slower than $r$ as $r\to 0$.
   
   Our first result  generalizes Theorem 1 in \cite{MM1}, which only applied  to ${\mathcal L}^o(x,\partial_x)$.
      
 \begin{theorem}\label{th:1}  
 For $n\geq 3$ and $p\in (1,\infty)$, assume $a_{ij}=a_{ji}$ satisfy \eqref{aij-delta} while the $b_k$ satisfy \eqref{est:b_k}.
 Then  for $\e>0$ sufficiently small, there is a solution of \eqref{LZ=0} in the form
 \begin{subequations}\label{Z(x)}
 \begin{equation}\label{Z(x)=}
 Z(x)=\int_{|x|}^\e s^{1-n}\,e^{I(s)}\,ds\,(1+\zeta(|x|)) + v(x) \quad\hbox{for}\ |x|<\e,
 \end{equation}
 where $I(r)$ is given by \eqref{def:I}, $\zeta(r)$ satisfies
 \begin{equation}\label{est:M_{2,p}(zeta)}
 M_{2,p}(\zeta,r)\leq c\,\max(\om(r),\sigma(r)),
 \end{equation}
 and $v(x)$ satisfies
  \begin{equation}\label{est:M(v)}
 M_{2,p}(v,r)\leq c\,r^{2-n}\,e^{I(r)}\,\om(r).
  \end{equation}
  \end{subequations}
  Moreover,
for any $u\in W^{2,p}_\loc(\overline{B_\e}\backslash\{0\})$
that is a strong solution of ${\mathcal L}(x,\partial_x)u=0$ in $\overline{B_\e}\backslash\{0\}$ subject to the growth condition
 \begin{equation}\label{u-growthcondition}
 M_{2,p}(u,r)\leq c\,r^{1-n+\e_0}\ \hbox{where}\ \e_0>0,
 \end{equation}
 there exist constants $c_0$, $c_1$ (depending on $u$) such that
  \begin{equation}\label{u-expansion}
  u(x)=c_0\,Z(x)+ c_1+w(x),
   \end{equation}
   where $M_{2,p}(w,r)\leq c\,r^{1-\e_1}$ for any $\e_1>0$.
\end{theorem}
\noindent
We will prove this result in Section \ref{sec:1}, but let us here observe that $I(r)$ satisfies $|I'(r)|\leq c\,r^{-1}\,\om(r)$, so we may integrate by parts and take $p>n$ to conclude
\begin{equation}\label{Z-asym}
Z(x)=\frac{|x|^{2-n}\,e^{I(|x|)}}{n-2}\left(1+\xi(x)\right),
\end{equation}
where $M_{1,\infty}(\xi,r)\leq c\,\max(\om(r),\sigma(r))$. 
 This shows that $I(r)$ controls how closely $Z$ adheres to the fundamental solution of the Laplacian.

Our second result generalizes Theorem 2 in \cite{MM1} and shows that the existence and finiteness of the limit $I(0)=\lim_{r\to 0}I(r)$ determines whether $Z$ solves \eqref{LZ=C0*delta} for some constant $C_0$.

 \begin{theorem}\label{th:2}   
 Under the assumptions of Theorem \ref{th:1} and $\e>0$ sufficiently small:
         \begin{itemize}
\item[(i)] If $I(0)=\lim_{r\to 0}I(r)$ exists and is finite, then
we can solve \eqref{LZ=C0*delta} in $B_\e$ with $C_0=|S^{n-1}|\,e^{I(0)}$.
\item[(ii)] If $I(r)\to -\infty$ as $r\to 0$, then solving  \eqref{LZ=C0*delta} in $B_\e$ yields $C_0=0$, and so $Z$ solves $-{\mathcal L}(x,\partial_x) Z(x)=0$ in $B_\e$, despite its singularity at $x=0$.
\end{itemize}
 \end{theorem}
 \noindent
 We will prove this result in Section \ref{sec:2}, but let us observe that in case (i) we have found a solution of
\eqref{-LF=delta} of the form
 \begin{equation}\label{F-asym}
 F(x)=\frac{|x|^{2-n}}{(n-2)|S^{n-1}|}\,(1+o(|x|)) \quad\hbox{as}\ |x|\to 0,
 \end{equation}
 so the comparison with the fundamental solution of the Laplacian is even more explicit.
On the other hand, if  $I(r)\to +\infty$ as $r\to 0$, then the singular solution $Z(x)$ grows more rapidly as $|x|\to 0$ than the fundamental solution for the Laplacian, and we are not able to solve \eqref{LZ=C0*delta}. 
 
 Let us consider a simple example to illustrate the effect of the drift term on the fundamental solution.
 
 \medskip\noindent
 {\bf Example 1.} 
 Let $a_{ij}=\delta_{ij}$ so that our operator \eqref{L-nondivergence} becomes
 \begin{equation}\label{L:Ex1}
 {\mathcal L}(x,\partial_x)=\Delta + b_k(x)\partial_x,
 \end{equation}
 and the quantity $I(r)$ in \eqref{def:I} reduces to
  \begin{equation}\label{I(r):Ex1}
  I(r)=\frac{1}{|S^{n-1}|}\int_{r<|x|<\e} \frac{\langle {\bf B}_z,z\rangle}{|z|^n}\,dz.
  \end{equation} 
  If $\om(r)$ satisfies the Dini condition, then from condition \eqref{est:b_k} we easily conclude that $I(0)$ exists and is finite; but this finite limit may exist without the Dini condition: e.g.\ we could take $b_k(x)=x_k\,\sin(|x|^{-1})\,\om(|x|)/|x|^2$. In any case, provided $I(0)$ is finite, we have a solution $F(x)$ of \eqref{LZ=C0*delta} that is comparable to the fundamental solution of the Laplacian.
    
In \cite{CZ},  Cranston and Zhao consider operators of the form ${\mathcal L}=\frac{1}{2}\Delta+b\cdot\nabla$ with vector field $b$. 
Assuming $U$ is a bounded Lipschitz domain and $b(x)$ satisfies the conditions
\begin{subequations}
\begin{equation}\label{CZ-cond}
\lim_{r\to 0}\,\sup_{x\in U}\int_{|x-y|<r}\frac{|b(y)|^2}{|x-y|^{n-2}}dy=0=\lim_{r\to 0}\,\sup_{x\in U}\int_{|x-y|<r}\frac{|b(y)|}{|x-y|^{n-1}}dy,
\end{equation}
they conclude that the Green's function $G(x,y)$ for ${\mathcal L}$ in $U$ exists and is comparable to the Green's function $G_0(x,y)$ for ${\mathcal L}_0=\frac{1}{2}\Delta$, i.e.\ 
\begin{equation}\label{G,G0-equiv}
c^{-1}\,G_0(x,y)\leq      G(x,y)\leq c\,G_0(x,y) \quad\hbox{for}\ x,y\in U, \ x\not=y.
\end{equation}
\end{subequations}
If $U$ contains the origin, then setting $F(x):=G(x,0)$ defines a solution of \eqref{-LF=delta}, so
 let us compare our hypotheses and conclusions with those of \cite{CZ}. If $b$ has a singularity at $x=0$ of the form
$|b(x)|=|x|^{-1}\,\om(|x|),$ the conditions in \eqref{CZ-cond} require
\[
\int_0^r \frac{\om^2(\rho)}{\rho}\,d\rho<\infty \quad\hbox{and}\quad
\int_0^r \frac{\om(\rho)}{\rho}\,d\rho<\infty.
\]
The first of these is the square-Dini condition \eqref{squareDini} that we have required, while the second is the Dini condition that we have {\it not} required: to conclude that the fundamental solution exists at $x=0$, we only require the function $I(r)$ given in \eqref{I(r):Ex1} to have a finite limit $I(0)=\lim_{r\to 0}I(r)$. Moreover, \eqref{Z-asym} is a sharper estimate than \eqref{G,G0-equiv}.

\medskip
In Example 1, the drift term plays the role of a perturbation which, if not too large, does not affect the existence of the fundamental solution.
We now consider  an example where singular drift can convert an operator ${\mathcal L}^o$ that does not have a fundamental solution to one for which a fundamental solution exists! 

  \medskip\noindent
 {\bf Example 2.} 
 Consider ${\mathcal L}^o=a_{ij}\partial_i\partial_j$ with coefficients
 \begin{equation}\label{GS-aij}
 a_{ij}(x)=\delta_{ij}+g(r)\frac{x_ix_j}{|x|^2},
 \end{equation}
 where $|g(r)|\leq\om(r)$ with $\om$ satisfying \eqref{squareDini}.
 Coefficients of the form \eqref{GS-aij} were first considered by Gilbarg \& Serrin \cite{GS}, and have proven useful in both nondivergence and divergence form equations (cf.\  \cite{MM1}, \cite{MM2},\cite{MM3}).
 We can use \eqref{def:I0} to calculate
  \begin{equation}
 I^o(r)=(1-n)\int_r^1 g(\rho)\frac{d\rho}{\rho}.
 \end{equation}
 Note that $I^o(0)\!=\!\lim_{r\to 0} I^o(r)$ exists and is finite when $\om$ satisfies the Dini condition, but this finite limit may exist without Dini continuity: e.g.\ $g(r)=\sin(r^{-1})\om(r)$.

Now let us assume $g(r)>0$ and $\om$ does not satisfy the Dini condition. Then $I^o(r)\to -\infty$ as $r\to 0$, so the fundamental solution for ${\mathcal L}^o$ does not exist at $y=0$.  
  However, if we add the first-order coefficients
  \begin{equation}\label{GS-b_k}
  b_k(x) =(n-1)\,\frac{x_k}{|x|^2} \,[g(|x|)+g^2(|x|)]
  \end{equation}
  to obtain ${\mathcal L}$ as in \eqref{L-nondivergence}, then we can use \eqref{def:I} to calculate
  \begin{equation}
  I(r)=(n-1)\int_r^1 \frac{g^2(\rho)}{\rho}\,d\rho.
  \end{equation}
  Since $g^2(r)\leq \om^2(r)$ and $\om$ satisfies  \eqref{squareDini}, we see that $I(r)$ has a finite limit as $r\to 0$, and so we can solve \eqref{LZ=C0*delta} to conclude the fundamental solution exists at $y=0$.

\medskip
 As previously stated, we will prove Theorem 1 in Section \ref{sec:1} and Theorem 2 in Section  \ref{sec:2}; but in those sections we will also state and prove Corollaries 1 and 2 respectively, which show how the formulas for the singular solution $Z$ change when $y\not= 0$ and $a_{ij}(y)\not=\delta_{ij}$.  To formulate these results, we need to generalize \eqref{aij-delta}: for
 a given $y\in U$,  choose $\e$ so that $0<\e<{\rm dist}(y,\partial U)$ and assume
 \begin{equation}\label{|A_x-A_y|<om}
\sup_{|x-y|=r} \|{\bf A}_x-{\bf A}_y\|\leq \om(r)\quad\hbox{for}\ 0<r<\e.
\end{equation}
We continue to assume the $b_k$ are only singular at $x=0$ and satisfy \eqref{est:b_k}.
We also generalize the function $I(r)$ in \eqref{def:I} as
\begin{subequations}\label{def:I_y}
 \begin{equation}\label{I_y=}
 I_y(r)=\frac{1}{|S^{n-1}|}\int_{r<|z-y|<\e} H(z,y) \,\frac{dz}{|z-y|^n}\, ,
 \end{equation}
 where the integrand $H(z,y)$ is
 \begin{equation}\label{H(x,y)=}
{\rm tr}({\bf A}_z{\bf A}_y^{-1})-n\frac{\langle {\bf A}_z{\bf A}_y^{-1/2}(z-y),{\bf A}_y^{-1/2}(z-y)\rangle}{|z-y|^2}+\langle {\bf B}_z{\bf A}_y^{-1/2},(z-y)\rangle,
 \end{equation}
 \end{subequations}
  with ${\bf A}_y^{-1}$ denoting the inverse matrix.
 If $y=0$ and ${\bf A}_0=I$, then \eqref{def:I_y} coincides with \eqref{def:I}; if $y\not=0$ and we also stipulate $0<\e<|y|$, then 
 ${\bf B}_z$ is bounded on $|z-y|<\e$, so the last term in \eqref{H(x,y)=} will not play a role in whether the limit $I_y(0)$ exists and is finite. Now we state the main result of this paper.
 \begin{theorem} 
Suppose ${\mathcal L}(x,\partial_x)$ as in \eqref{L-nondivergence} is an elliptic operator in a bounded open set $U\subset \RR^n$, $n\geq 3$, where the coefficients $a_{ij}=a_{ji}$ are continuous functions with modulus of continuity $\om(r)$ satisfying \eqref{squareDini}. 
Suppose $U$ contains the origin and the $b_k$ satisfy  \eqref{est:b_k} but otherwise are bounded in $U$.
For each $y\in U$ assume that the limit $I_y(0)=\lim_{r\to 0} I_y(r)$ exists and is finite. Then
${\mathcal L}(x,\partial_x)$ has a fundamental solution $F(x,y)$ in $U$ and it has the asymptotic behavior
 \begin{equation}\label{F(x,y)=}
 F(x,y)=\frac{\langle {\bf A}_y^{-1}(x-y),(x-y)\rangle^{\frac{2-n}{2}}}{(n-2)|S^{n-1}|\sqrt{\det{\bf A}_y}}(1+o(1)))\quad\hbox{as}\ x\to y.
 \end{equation}
\end{theorem}
\noindent
This is proved in Section 3. Note that the leading asymptotic in \eqref{F(x,y)=} is familiar as the ``Levi function'' that occurs in the classical case (cf.\ \cite{M}).
As a consequence of Theorem 3, we see that the singular drift may affect the existence of the fundamental solution but does not play a role in its asymptotic behavior as $x\to y$.

Many of the  arguments in this paper  also appear in \cite{MM1}, but we have repeated them here for the convenience of the reader.
From a more general perspective, the asymptotic analysis used here is related to that developed in \cite{KM}.

Let us compare our results (in \cite{MM1} and this paper) with recent work estimating the singularity of the Green's function for ${\mathcal L}^o(x,\partial_x)=a_{ij}(x)\partial_i\partial_j$ when the coefficients $a_{ij}$ satisfy the Dini mean oscillation (DMO) condition: cf.\ \cite{HK}, \cite{KL}, \cite{DKL}. If $U$ is a bounded $C^{1,1}$-domain and the $a_{ij}=a_{ji}$ are continuous functions on $\overline U$, then it is well-known that the Green's function $G(x,y)$ exists. If the $a_{ij}$ satisfy the DMO condition in $U$, then \cite{DKL} shows that, for any $x_0\in U$, the following limit holds:
\begin{equation}\label{DMO-asym}
\lim_{x\to x_0} |x-x_0|^{2-n} |G(x,x_0)-G_{x_0}(x,x_0)|=0,
\end{equation}
where $G_{x_0}$ denotes that Green's function for the constant coefficient operator ${\mathcal L}^o(x_0,\partial_x)$ in $U$.
Since the Green's function is a particular fundamental solution and since \eqref{DMO-asym} is equivalent to \eqref{F(x,y)=}, it is natural to compare the hypotheses of the two results. As shown in \cite{DK}, there are coefficients that are DMO but do not satisfy our square-Dini condition \eqref{squareDini}. On the other hand, there are coefficients of the form \eqref{GS-aij} which satisfy  \eqref{squareDini} but are not DMO. In fact, as shown in \cite{MM3}, with
\begin{equation}
g(r)=\sin(|\log r|)\,|\log r|^{-\gamma},
\end{equation}
the coefficients in  \eqref{GS-aij}  are DMO only for $\gamma>1$, but they satisfy  \eqref{squareDini} for $\gamma>1/2$ and the limit $I(0)$ exists and is finite for all $\gamma>0$. This example with $1/2<\delta\leq 1$ shows that the results of  \cite{DKL} do not cover our results 
for ${\mathcal L}^o(x,\partial_x)$, let alone the operator \eqref{L-nondivergence} with singular drift.
 
 
\section{Construction of the singular solution $Z$}\label{sec:1}

In this section, we will not only prove Theorem 1, but we will state and prove Corollary 1, which shows how the formulas change when we no longer assume $y=0$ and $a_{ij}(y)=\delta_{ij}$.  Instead of constructing $Z(x)$ in a small ball, we replace the condition that $\om(r)$ satisfies \eqref{squareDini} with
\begin{equation}\label{est:sigma(1)}
\sigma(1)=\int_0^1 \om^2(\rho) \frac{d\rho}{\rho}<\mu^2
\end{equation}
where $\mu>0$ is sufficiently small, and then show existence in the unit ball $B_1$. In fact, with $\kappa\in (0,1)$ as in 
 \eqref{om-monotone}, this also implies
\begin{equation}\label{est:om(r)}
\om(r)<c_\kappa\,\mu\quad\hbox{for}\ 0<r\leq 1,
\end{equation}
since
\[
\mu^2>\int_{0}^r \om^2(\rho)  \frac{d\rho}{\rho}\geq \om^2(r) r^{-2+2\kappa}\int_{0}^r \rho^{1-2\kappa}d\rho=\frac{\om^2(r)}{2(1-\kappa)}.
\]
Moreover, it will be useful to consider ${\mathcal L}$ on all of $\RR^n$, so we assume
\begin{equation}\label{a_ij=delta_ij}
a_{ij}(x)=\delta_{ij} \ \hbox{and} \ b_k(x)=0 \ \hbox{for}\ |x|>1,
\end{equation}
and construct a solution of ${\mathcal L}Z=0$ on $\RR^n\backslash\{0\}$.
  
  \medskip\noindent
  {\bf Proof of Theorem 1.} 
As in \cite{MM1}, we use spherical means: for a function $f(x)$ we denote its mean value over the sphere $|x|=r$ by $\overline f(r)$:
\[
\overline{f}(r)= \meanint_{S^{n-1}} f(r\theta)\,d\theta,
\]
where $S^{n-1}$ is the unit sphere, the slashed integral denotes mean value, $r=|x|$, $\theta=x/|x|\in S^{n-1}$, and $d\theta$ denotes standard surface measure on $S^{n-1}$.   
Let us write
  \begin{subequations}
  \begin{equation}\label{Z=h+v}
  Z(x)=h(|x|)+v(x), \quad\hbox{where} \ h(r):=\overline{Z}(r),
  \end{equation}
  so that
  \begin{equation}
  \overline{v}(r)=0.
  \end{equation}
  \end{subequations}
 If we take the spherical mean of the equation ${\mathcal L}(h+v)=0$, we obtain
 \begin{subequations}
\begin{equation}\label{ode-h-1}
 \alpha(r)\,h''+\left[\frac{\alpha_n(r)-\alpha(r)+\beta(r)}{r}\right]\,h'+\overline{{\mathcal L}\,v}(r)=0,
\end{equation}
where
       \begin{equation}\label{def:alpha}
     \alpha(r):=\meanint_{S^{n-1}} a_{ij}(r\theta)\theta_i\theta_j\,d\theta, \qquad    \alpha_n(r):=\meanint_{S^{n-1}} a_{ii}(r\theta)\,d\theta, \\
   \end{equation}
and
\begin{equation}\label{def:beta0}
   \beta(r):= r\,\meanint_{S^{n-1}} b_k(r\theta)\theta_k\,d\theta.
\end{equation}
\end{subequations}
In terms of these, instead of \eqref{def:I} we can write
\begin{equation}\label{I(r)=}
I(r)=\int_r^1 (\alpha_n(s)-n\,\alpha(s)+\beta(s))\frac{ds}{s}.
\end{equation}
Using \eqref{aij-delta}, we see that
\begin{subequations}
        \begin{equation}\label{def:alpha-limits}
     \alpha(r)=1+O(\om(r)) \quad \hbox{and} \quad    \alpha_n(r)=n+O(\om(r))\quad\hbox{as}\ r\to 0, \\
   \end{equation}
   and using \eqref{est:b_k} we see that
   \begin{equation}\label{est:beta}
  |\beta(r)| \leq \om(r).
\end{equation}
\end{subequations}
Hence the integrand in \eqref{I(r)=} is bounded by $c\,\om(s)\,s^{-1}$. Since $\om(r)$ need not satisfy the Dini condition, we do not know whether $I(r)$ has a finite limit as $r\to 0$; this is, of course, the significance of Theorem 2.

Turning to the $\overline{{\mathcal L}\,v}(r)$ term, if we write $a_{ij}\partial_i\partial_j v=\widetilde a_{ij}\partial_i\partial_j v+\Lap v$ where $\widetilde a_{ij}:=(a_{ij}-\delta_{ij})$ and use $\overline{\Lap v}=\Lap \overline{v}=0$, we see that
\[
\overline{{\mathcal L}v}(r)=\overline{\widetilde a_{ij}\partial_i\partial_j v}(r)+\overline{ b_k\partial_k v}(r).
\]
Notice that $|\overline{{\mathcal L}v}(r)|\leq c\,\om(r)\,(|\overline{D^2v}|+r^{-1}|\overline{Dv}|)$ for $0<r<1$
and $\overline{{\mathcal L}v}(r)=0$ for $r>1$.
Since $\alpha(r)\to 1$ as $r\to 0$ and $\alpha(r)=1$ for $r>1$, we may assume $\alpha(r)\geq \delta>0$ for $0<r<\infty$. Hence we may divide \eqref{ode-h-1} by $\alpha(r)$ and replace $h'$ by $g$ to obtain
      \begin{subequations}
 \begin{equation}\label{ode-h-3}
g'+\left[\frac{n-1+R(r))}{r}\right]\,g= B[v](r),
 \end{equation}
where 
 \begin{equation}
 R(r)=\frac{\alpha_n(r)+\beta(r)}{\alpha(r)}-n 
 \end{equation}
 satisfies $|R(r)|\leq c\,\om(r)$ as $r\to 0$ and $R(r)=0$ for $r>1$,
and $B[v](r)$ satisfies 
  \begin{equation}\label{est:B[v]}
|B[v](r)|\leq c\,\om(r) \left(|\overline{D^2 v}(r)|+r^{-1}|\overline{Dv}(r)|\right) \quad\hbox{for}\ 0<r<1
 \end{equation}
 and $B[v](r)=0$ for $r>1$. Moreover, the monotonicity of $\om(r)$ and \eqref{om-monotone} imply
 \begin{equation}\label{<om(r)<}
 \sup_{r<\rho<2r}\om(\rho)\leq c\,\om(r),
 \end{equation}
 so we consequently obtain
\begin{equation}\label{est:M_p(B[v])}
r^2 M_p(B[v],r)\leq c\,\om(r)\, M_{2,p}(v,r).
\end{equation}
   \end{subequations} 
   
   Solving \eqref{ode-h-3} involves the integrating factor $r^{n-1}E_-(r)$, where we introduce
   \begin{equation}\label{def:E_pm}
   E_\pm(r)=\exp\left[\pm\int_r^\infty R(s)\,\frac{ds}{s}\right]=\exp\left[\pm\int_r^1 R(s)\,\frac{ds}{s}\right]=\frac{1}{E_\mp(r)}.
   \end{equation}
    \begin{subequations}
   Notice that $E_-(r)E_+(\rho)=\exp(\int_\rho^r R(s)\,s^{-1}\,ds)$ so by \eqref{est:om(r)} we have
   \begin{equation}\label{est:<exp(int_rho^r)}
   \left(\frac{\rho}{r}\right)^{c_\kappa\,\mu}\leq \exp\left( \pm \int_\rho^r R(s)\frac{ds}{s}\right)\leq \left(\frac{r}{\rho}\right)^{c_\kappa\,\mu}
   \quad\hbox{for}\ 0<\rho<r\leq 1.
   \end{equation}
  In particular, we have
   \begin{equation}\label{<E_-(r)<}
   c_1\,E_\pm(r)\leq E_\pm(\rho)\leq c_2\,E_\pm(r)\quad\hbox{for}\ r<\rho<2r,
   \end{equation}
   and for any $f\in L^p_{\loc}(\RR^n\backslash\{0\})$ and $\nu\in\RR$ we have
    \begin{equation}
    M_p(|x|^\nu E_\pm(|x|)\,f(x),r)\leq c\,r^\nu E_\pm(r)\,M_p(f,r).
     \end{equation}
   \end{subequations}
   While $E_+(r)$ is used to solve  \eqref{ode-h-3}, we observe that it is equivalent to $e^{I(r)}$. In fact,
   \begin{subequations}
   \begin{equation}\label{E+=A*E^(I)}
   E_+(r)=A\,e^{I(r)}(1+\tau(r))
   \end{equation}
   where $|R(s)(1-\alpha(s))|\leq c\,\om^2(s)$ implies
     \begin{equation}
     A=\exp\left[\int_0^1 R(s)(1-\alpha(s))s^{-1}\,ds\right] \quad\hbox{is finite and positive,}
     \end{equation}
   and
   \begin{equation}\label{tau=}
   \tau(r)=\exp\left[ -\int_0^r R(s)[1-\alpha(s)]\,\frac{ds}{s}\right]-1 \quad\hbox{satisfies} \ |\tau(r)|\leq c\,\sigma(r).
   \end{equation}
    \end{subequations}
 Hence, for some constants $c_1,c_2$ we have
    \begin{equation}\label{E_+=exp(I)}
    c_1\,E_+(r)\leq e^{I(r)}\leq c_2\,E_+(r).
     \end{equation}
     
We consider \eqref{ode-h-3} as an ODE that depends on $v\in Y$, where $Y$ is the Banach space of functions $v\in W^{2,p}_{\loc}(\RR^n\backslash\{0\})$ for which 
 \begin{equation}\label{def:Y}
 \|v\|_Y:=\sup_{0<r<1}\frac{M_{2,p}(v,r)\,r^{n-2}}{\om(r)\,e^{I(r)}} + \sup_{r>1}\frac{M_{2,p}(v,r)\,r^{n-1}}{\mu} <\infty.
 \end{equation}     
If we let $\phi(r)=r^{n-1}E_-(r)\,g(r)$, then \eqref{ode-h-3} implies
$
\phi'(r)=r^{n-1} E_-(r)\,B[v](r).
$
This can be integrated to find
\begin{equation}\label{phi(r)=}
\phi(r)=\phi(0)+\int_0^r \rho^{n-1} E_-(\rho)\,B[v](\rho)\,d\rho,
\end{equation}
where $\phi(0)$ is an arbitrary constant. Of course, for \eqref{phi(r)=} to be valid we need to know that $r^{n-1}E_-(r)\,B[v](r)$ is integrable at $r=0$. In fact, we will show below that for $v\in Y$ we have
  \begin{equation}\label{est:integral}
  \int_0^r \rho^{n-1}\,E_-(\rho)\,|B[v](\rho)|\,d\rho\leq c\,\mu^2\quad\hbox{for all} \ r>0,
  \end{equation}
  where $c$ may be taken uniformly for all $v\in Y$ with $\|v|_Y\leq 1$.
Hence \eqref{phi(r)=} is valid and we conclude
  \begin{subequations}
  \begin{equation}\label{h'(r)=}
  g(r)=h'(r)=r^{1-n}\,E_+(r)\left[\phi(0)+\int_0^r\rho^{n-1}\,E_-(\rho)\,B[v](\rho)\,d\rho\right]
  \end{equation}
  and
  \begin{equation}\label{h''(r)=}
  h''(r)=\frac{1-n-R(r)}{r^n}\,E_+(r)\left[\phi(0)+\int_0^r\rho^{n-1}E_-(\rho) B[v](\rho)\,d\rho\right]+B[v](r).
  \end{equation}
  \end{subequations}
  
  To verify \eqref{est:integral}, we observe that $v\in Y$ implies $M_{2,p}(v,r)\leq c\,\om(r)\,r^{2-n}\,E_+(r)$, so we can use \eqref{est:B[v]}, \eqref{<om(r)<},  \eqref{<E_-(r)<}, and H\"older's inequality to conclude
  \[
  \begin{aligned}
  \int_r^{2r} \rho^{n-1}\,E_-(\rho)\,|B[v](\rho)|\,d\rho 
  & \leq  c\,E_-(r)\,\om(r)\int_r^{2r}\rho^{n-1}\left( |\overline{D^2 v}(\rho)|+\rho^{-1}|\overline{Dv}(\rho)|\right)d\rho \\
  & \leq  c\,E_-(r)\om(r)\int_{r<|x|<2r}\left(|D^2v(x)|+|x|^{-1}|Dv(x)|\right) dx\\
  & \leq c\,E_-(r)\,\om(r)\,r^{n-2}\,M_{2,p}(v,r) \leq c\,\om^2(r).
  \end{aligned}
  \]
  Now if we write
  \[
  \int_0^r \rho^{n-1}\,E_-(\rho)\,|B[v](\rho)|\,d\rho=\sum_{j=0}^\infty \int_{r/2^{j+1}}^{r/2^j} \rho^{n-1}\,E_-(\rho)\,|B[v](\rho)|\,d\rho,
  \]
  then we obtain 
  \begin{equation}\label{est:int-B[v]}
  \begin{aligned}
  \int_0^r \rho^{n-1}\,E_-(\rho)\,|B[v](\rho)|\,d\rho & \leq c\sum_{j=0}^\infty \om^2\left(\frac{r}{2^{j+1}}\right) \\
 & \leq c\int_0^r\om^2(\rho)\frac{d\rho}{\rho}=c\,\sigma(r)<c\,\mu^2.
  \end{aligned}
  \end{equation}
  This confirms  \eqref{est:integral} with $c$ uniform for $\|v\|_Y\leq 1$.
 
 We also have a PDE for $v$ that depends upon $h$. From ${\mathcal L}Z-\overline{{\mathcal L}Z}=0$ we find:
  \begin{equation}\label{pde-v}
  \begin{aligned}
  -\Lap v= \, & \widetilde{a}_{ij}\partial_i\partial_j h-\overline{ \widetilde{a}_{ij}\partial_i\partial_j h}+\widetilde{a}_{ij}\partial_i\partial_j v-\overline{ \widetilde{a}_{ij}\partial_i\partial_j v}\\
  &+b_k\partial_k h - \overline{b_k\partial_k h}+
  +b_k\partial_k v -\overline{ b_k\partial_kv }.
  \end{aligned}
   \end{equation}
For a given $v\in Y$, we solve \eqref{ode-h-3} for $h$ and use \eqref{h'(r)=} and  \eqref{h''(r)=} to write
\[
b_k\partial_k h=r^{-n}E_+(r)\left[\phi(0)+\int_0^r\rho^{n-1}E_-(\rho)B[v](\rho)d\rho\right]\psi_1
\]
and
\[
\widetilde a_{ij}\partial_i\partial_j h=r^{-n}E_+(r)\left[\phi(0)+\int_0^r \rho^{n-1}E_-(\rho) B[v](\rho)\,d\rho\right]\,\psi_2 +B[v]\,\widetilde a_{ij}\theta_i\theta_j,
\]
where
\begin{equation}
\psi_1(r\theta)=r\,b_k(r\theta)\theta_k \quad\hbox{and}\quad
\psi_2(r\theta)=\widetilde a_{ii}(r\theta)-(n+R(r))\widetilde a_{ij}(r\theta)\theta_i\theta_j
\end{equation}
also satisfiy $|\psi_i(r\theta)|\leq c\,\om(r)$ for $i=1,2$. Plugging this into \eqref{pde-v}, we obtain
 an equation of the form  $-\Lap v= F[v]$. We want to apply $K$, convolution by the fundamental solution of the Laplacian, to solve this, but there could be a problem:  $F[v]$ may not be integrable at $x=0$. However, from \eqref{pde-v} we see that $\overline{F[v]}=0$ and multiplying $F[v](x)$ by $|x|$ makes it integrable, so we can use Proposition \ref{pr:convolution} in Appendix A.

Applying $K$ to both sides of  \eqref{pde-v}, we obtain an equation for $v$ alone:
\begin{equation}\label{OperatorEqn}
v+S_1v+S_2v=\phi(0)w,
\end{equation}
where
\[
w(x)=K_{y\to x}\left( |y|^{-n} E_+(|y|)\,\left[ \psi(y)-\overline{\psi}(|y|)\right]\right)
\]
\[
S_1v\!=\!-K_{y\to x}\left(\! |y|^{-n}E_+(|y|)\int_0^{|y|}\rho^{n-1}E_-(\rho)B[v](\rho)d\rho\left[\psi(y)-\overline{\psi}(|y|)\right]\!\right)
\]
\[
S_2v=-K_{y\to x}\left( B[v](|y|)\left[\widetilde a_{ij}\theta_i\theta_j-\overline{\widetilde a_{ij}\theta_i\theta_j}\right]+\widetilde a_{ij}\theta_i\theta_j-\overline{\widetilde a_{ij}\theta_i\theta_j}\right),
\]
and we have let $\psi:=\psi_1+\psi_2$.
To find $v\in Y$ satisfying \eqref{OperatorEqn} we need to show $w\in Y$ and $S_i:Y\to Y$ has small operator norm for $i=1,2$.

To show $w\in Y$, we must estimate $M_{2,p}(w,r)$ for $0<r<1$ and for $r>1$. We apply Proposition 1 in Appendix \ref{App1} to $f(x)=|x|^{-n} E_+(|x|)(\psi(x)-\overline{\psi}(|x|)$, which vanishes for $|x|>1$, to conclude
\begin{equation}\label{est:M_2(w)}
M_{2,p}(w,r)\!\leq \!c\left(\! r^{1-n}\!\int_0^r E_+(\rho)\om(\rho)\,d\rho\!+\!r\!\int_r^1 E_+(\rho)\om(\rho)\rho^{-n}\,d\rho\!\right) \ \hbox{for}\ 0<r<1.
\end{equation}
We can use \eqref{est:<exp(int_rho^r)} and the monotonicity of $\om(r)$ to estimate
\begin{subequations}
\begin{equation}\label{est:int_0^r-E+*om}
\int_0^r E_+(\rho)\,\om(\rho) \,d\rho\leq E_+(r)\,\om(r) \,r^{c\mu}\int_0^r\rho^{-c\mu}d\rho =
c\, r\, E_+(r)\,\om(r)
\end{equation}
and similarly
\begin{equation}\label{est:int_r^1-E+*om}
\int_r^1 E_+(\rho)\,\om(\rho)\,\rho^{-n}\,d\rho\leq c\,r^{1-n}\,E_+(r)\,\om(r).
\end{equation}
\end{subequations}
Using these in \eqref{est:M_2(w)}, we obtain
$M_{2,p}(w,r)\,r^{n-2}\leq c\,\om(r)\,E_+(r)\quad\hbox{for}\ 0<r<1,$
and we can invoke \eqref{E_+=exp(I)} to replace $E_+(r)$ by $e^{I(r)}$ as required in the norm for $Y$. 
Meanwhile, for $r>1$ we
use \eqref{est:om(r)} and $E_+(\rho)\leq \rho^{-c\,\mu}$ for $0<\rho<1$ to conclude
\[
M_{2,p}(w,r)\leq c\,r^{1-n}\int_0^1 E_+(\rho)\,\om(\rho)\,d\rho\leq c\,\mu\,r^{1-n}.
\]
Hence $M_{2,p}(w,r)\,r^{n-1}\leq c\,\mu \quad\hbox{for} \ r>1.$
These estimates confirm that $w\in Y$.

Next let us show that $S_1:Y\to Y$ with small operator norm. We assume $\|v\|_Y\leq 1$ and we want to estimate $M_{2,p}(S_1v,r)$ separately for $0<r<1$ and $r>1$. Using \eqref{est:integral} we see that the function
\[
f_1(y)=|y|^{-n} E_+(|y|)\int_0^{|y|}\rho^{n-1} E_-(\rho)\,B[v](\rho)\,d\rho\left(\psi(y)-\overline{\psi}(|y|)\right)
\]
satisfies $M_p(f_1,r)\leq c\,\mu^2\,E_+(r)\,\om(r)\,r^{-n}$ for $0<r<1$ and $M_p(f_1,r)=0$ for $r>1$.
Since $S_1v=-Kf_1$, we can apply Proposition \ref{pr:convolution} in Appendix 1 to obtain
\[
M_{2,p}(S_1v,r)\leq c\,\mu^2\left(r^{1-n}\int_0^r E_+(\rho)\om(\rho)\,d\rho+r\int_r^1E_+(\rho)\om(\rho)\rho^{-n}\,d\rho\right).
\]
Using \eqref{est:int_0^r-E+*om} and \eqref{est:int_r^1-E+*om}, we conclude $M_{2,p}(S_1 v,r)\,r^{n-2}\leq c\,\mu^2\,\om(r)\,E_+(r)$ for $0<r<1$. On the other hand, for $r>1$, Proposition 1 in Appendix 1 implies
\[
M_{2,p}(S_1 v,r)\leq c\,\mu^2\,r^{1-n}\int_0^1 E_+(\rho)\,\om(\rho)\,d\rho\leq c\,\mu^3\,r^{1-n},
\]
so $M_{2,p}(S_1v,r)\,r^{n-1}\leq c\,\mu^3$ for $r>1$. Combining these estimates, we see that $S_1:Y\to Y$ has small operator norm.

Finally, we show that $S_2:Y\to Y$ with small operator norm. Again we assume $\|v\|_Y\leq 1$ and estimate $M_{2,p}(S_2v,r)$ separately for $0<r<1$ and $r>1$. Notice that the function
\[
f_2=B[v]\,(\widetilde a_{ij}\theta_i\theta_j-\overline{\widetilde a_{ij}\theta_i\theta_j})
\]
satisfies
\[
M_p(f_2,r)\leq \om(r) M_p(B[v],r)\leq c\,\om^3(r)\,E_+(r)\,r^{-n} \quad\hbox{for}\ 0<r<1,
\]
where $c$ is independent of $v$, and $M_p(f_2,r)=0$ for $r>1$. Similarly, the function
\[
f_3=\widetilde a_{ij}\partial_i\partial_j v-\overline{\widetilde a_{ij}\partial_i\partial_j v}
\]
satisfies
\[
M_p(f_3,r)\leq \om(r)\,M_p(D^2v,r)\leq \om^2(r)\,E_+(r)r^{-n}\quad\hbox{for}\ 0<r<1,
\]
and $M_p(f_3,r)=0$ for $r>1$. For $0<r<1$, we apply Proposition 1 in Appendix 1 to $S_2v=-K(f_2+f_3)$ to conclude
\[
M_{2,p}(S_2 v,r)\leq c\left( r^{1-n}\int_0^r \om^2(\rho)\,E_+(\rho)\,d\rho+r\int_r^1 \om^2(\rho)\,E_+(\rho)\,\rho^{-n}\,d\rho\right).
\]
Using \eqref{est:om(r)}, \eqref{est:int_0^r-E+*om}, and \eqref{est:int_r^1-E+*om}, we conclude that $M_{2,p}(S_2 v,r)\,r^{n-2}\leq c\,\mu\,\om(r)\,E_+(r)$ for $0<r<1$. Meanwhile, for $r>1$, we use \eqref{est:om(r)} and \eqref{est:<exp(int_rho^r)} to estimate
\[
M_{2,p}(S_2v,r)\leq c\,r^{1-n}\int_0^1\om^2(\rho)E_+(\rho)\,d\rho\leq c\,\mu^2\,r^{1-n}\int_0^1\rho^{-c\,\mu}\,d\rho\leq c\,\mu^2\,r^{1-n}.
\]
Consequently, $M_{2,p}(S_2v,r)\,r^{n-1}/\mu\leq c\,\mu$. These estimates show that $S_2:Y\to Y$ has small operator norm.

Since  $S_1+S_2$ has small operator norm on $Y$, we conclude that \eqref{OperatorEqn} has a unique solution $v\in Y$, depending on the choice of the constant $c_*=\phi(0)$. But once $c_*$ and $v$ are known, we find $g(r)$ from \eqref{h'(r)=} and integrate to find $h(r)$:
\begin{equation}\label{h(r)=}
h(r)=\int_r^\infty s^{1-n}E_+(s)\left[c_*+\int_0^s\rho^{n-1}E_-(\rho)B[v](\rho)\,d\rho\right]ds+c_2
\end{equation}
where $c_2$ is an arbitrary constant. To show that the solution $Z(x)=h(|x|)+v(x)$ is of the form \eqref{Z(x)}, we choose $c_*$ to enable us to replace $E_+(s)$ by $e^{I(s)}$; recalling \eqref{E+=A*E^(I)} we see that we should take $c_*=A^{-1}$ and write $h(r)=h_0(r)+h_1(r)+c$ where 
\begin{subequations}\label{h0,h1}
\begin{equation}
h_0(r)=\int_r^1 s^{1-n}\,e^{I(s)}\,ds
\end{equation}
and (recalling \eqref{tau=})
\begin{equation}
h_1(r)=\int_r^1 s^{1-n}e^{I(s)}\tau(s)\,ds+\int_r^1 s^{1-n}E_+(s)\int_0^s\rho^{n-1}E_-(\rho)B[v](\rho)\,d\rho\,ds.
\end{equation}
\end{subequations}
Integrating by parts and using $|I'(s)|\leq c\,s^{-1}\om(s)$ we can show
\[
\left|h_0(r)-\frac{r^{2-n}}{n-2}e^{I(r)}\right|\leq c\,r^{2-n}\,e^{I(r)}\,\om(r).
\]
Using $|\tau(s)|\leq c\,\si(s)$ and \eqref{est:int-B[v]}, we can also estimate
\[
|h_1(r)|\leq c\,r^{2-n}\,e^{I(r)}\max(\om(r),\si(r))\quad\hbox{for}\ 0<r<1.
\]
If we define
\begin{equation}\label{def:zeta}
\zeta(r)=\frac{h_1(r)}{h_0(r)}\quad\hbox{for}\ 0<r<1,
\end{equation}
then we can easily estimate $|\zeta(r)|,|r\zeta'(r)|\leq c\,\max(\om(r),\si(r))$. To estimate $\zeta''$, let us write $h_0\zeta''=h_1''-h_0''\zeta-2h_0'\zeta'$ where 
\begin{subequations}
\begin{equation}\label{h0''=}
h_0''(r)=(n-1)r^{-n}e^{I(r)}-r^{1-n}e^{I(r)}I'(r)
\end{equation}
and
\begin{equation}\label{h1''=}
\begin{aligned}
h_1''(r)&=r^{-n}e^{I(r)}[(n-1)\tau(r)-rI'(r)\tau(r)-r\tau'(r)]\\
&+r^{-n}E_+(r)(n-1+R(r))\int_0^r\rho^{n-1}\,R(\rho)\,B[v](\rho)\,d\rho -B[v](r).
\end{aligned}
\end{equation}
\end{subequations}
We can estimate $h_0''$ and $h_0'$ pointwise, but $h_1''$ needs to be estimated in $M_p$. However, using \eqref{est:M_p(B[v])} and $v\in Y$, we may conclude $M_p(r^2\zeta'',r)\leq c\,\max(\om(r),\si(r))$. Combining with the estimates of the lower-order derivatives, we have shown \eqref{est:M_{2,p}(zeta)}. Since \eqref{est:M(v)} follows from $v\in Y$, we have proved the first part of Theorem 1.

Let us now turn to the second part of Theorem 1. Suppose $u\in W^{2,p}_{\loc}(\overline{B_1}\backslash\{0\})$ satisfies ${\mathcal L}u=0$ in $\overline{B_1}\backslash\{0\}$ and $M_{2,p}(u,r)\leq c\,r^{1-n+\e_0}$ for some $\e_0>0$; we want to show $u$ is of the form \eqref{u-expansion}. We will use properties of the bounded linear map 
\begin{equation}\label{L:W->L}
{\mathcal L}: W^{2,p}_{\delta_0,\delta_1}(\RR^n_o)\to L^{p}_{\delta_0+2,\delta_1+2}(\RR^n_o),
\end{equation}
where $\RR^n_o=\RR^n\backslash\{0\}$ and $W^{2,p}_{\delta_0,\delta_1}(\RR^n_o), L^{p}_{\delta_0+2,\delta_1+2}(\RR^n_o)$ are the weighted Sobolev spaces that are defined in Appendix \ref{App2}. Since we are interested in the behavior of functions at the origin, we fix $\delta_1\in (-n/p,-2+n/p')$ and allow $\delta_0$ to vary.
Note that \eqref{L:W->L} is a perturbation of \eqref{Lap:W->L}, 
 and the norm of the difference ${\mathcal L}-\Lap$ depends on the magnitude of 
\[
\sup_{|x|<1} \left(|a_{ij}(x)-\delta_{ij}| + |x|\,|b_k(x)| \right).
\]
Thus, provided we take $\mu$ in \eqref{est:sigma(1)} sufficiently small, we can arrange that \eqref{L:W->L} and \eqref{Lap:W->L} are not only Fredholm for exactly the same values of $\delta_0$ and $\delta_1$, but the nullity and deficiency of the two maps agree. Hence, with 
$\delta_1$ satisfying \eqref{delta1-range}, we find 
\begin{itemize}
\item[(i)] \eqref{L:W->L} is an isomorphism for $\delta_0\in (-n/p,-2+n/p')$;
\item[(ii)] \eqref{L:W->L} is surjective with nullity $1$ for $\delta_0\in (-2+n/p',-1+n/p')$;
\item[(iii)] \eqref{L:W->L} is injective with deficiency $1$ for $\delta_0\in (-n/p-1,-n/p)$.
\end{itemize}

Introduce a cutoff function $\chi\in C_0^\infty(B_1)$ with $\chi=1$ on $B_{1/2}$, then $M_{2,p}(\chi u,r)\leq c\,r^{1-n+\e_0}$ implies $\chi u\in W^{2,p}_{\delta_0}(B_1)$ provided $\delta_0>-1-\e_0+n/p'$. Let us choose $\delta^+_0\in (-1-\e_0+n/p',-1+n/p')$ and let $f={\mathcal L}(\chi u)$. Since $f=0$ for $|x|<1/2$ and for $|x|>1$, we have $f\in L^p_{\delta_0+2,\delta_1+2}(\RR^n_o)$ for all $\delta_0$ so let us choose $\delta_0\in (-n/p,-2+n/p')$. By (i) we can find $v={\mathcal L}^{-1}f\in W_{\delta_0,\delta_1}^{2,p}(\RR^n_o)$. Then $\chi u-v\in W^{2,p}_{\delta_0^+,\delta_1}(\RR^n_o)$ satisfies ${\mathcal L}(\chi u-v)=0$. By (ii), ${\mathcal L}:W^{2,p}_{\delta_0^+,\delta_1}(\RR^n_o)\to L^p_{\delta_0^++2,\delta_1+2}(\RR^n_o)$ has nullity 1. Since ${\mathcal L}Z=0$, the nullspace must be spanned by $Z$ and so $\chi u-v=c_0\,Z$ for some constant $c_0$. 

It only remains to show that $v=c_1+w$ where $M_{2,p}(w,r)\leq c\,r^{1-\e_1}$ for any $\e_1>0$.
Let us pick $\delta_0^-\in (-1-n/p,-n/p)$ so that by (iii), the map
\begin{equation}\label{L:W->L,codim}
{\mathcal L}: W^{2,p}_{\delta^-_0,\delta_1}(\RR^n_o)\to L^{p}_{\delta^-_0+2,\delta_1+2}(\RR^n_o),
\end{equation}
is injective with deficiency 1. Let $\zeta$ be a linear functional on $L^{p}_{\delta^-_0+2,\delta_1+2}(\RR^n_o)$ that vanishes on the image of 
\eqref{L:W->L,codim}. Note that ${\mathcal L}\chi=0$ for $|x|<1/2$ and for $|x|>1$, so ${\mathcal L}\chi\in L^p_{\delta_0^-+2,\delta_1+2}(\RR^n_o)$. But $\chi\not\in W^{2,p}_{\delta_0^-,\delta_1}(\RR^n_o)$ since $\delta_0^-<-n/p$, so ${\mathcal L}\chi$ is not in the image of \eqref{L:W->L,codim}, and hence $\zeta[{\mathcal L}(\chi)]\not=0$. This enables us to find $c_1$ so that
\[
\zeta[{\mathcal L}(\chi v)]=c_1\zeta[{\mathcal L}\chi].
\]
But this means $\zeta[{\mathcal L}(\chi v-c_1\chi)]=0$, i.e.\ ${\mathcal L}(\chi v-c_1\chi)$ is in the image of \eqref{L:W->L,codim}, so ${\mathcal L}w={\mathcal L}(\chi v-c_1\chi)$ for some $w\in W^{2,p}_{\delta^-_0,\delta_1}(\RR^n_o)\subset W^{2,p}_{\delta_0,\delta_1}(\RR^n_o)$. Since $\delta_0\in (-n/p,-2+n/p')$, by the isomorphism (i) we have $w=\chi v-c_1\chi$.   In other words, for $|x|<1/2$ we have 
$v=c_1+w$, where $w\in W^{2,p}_{\delta^-_0,\delta_1}(\RR^n_o)$. For any $\e_1\in (0,1)$ we can let $\delta_0^-=\e_1-1-n/p$ and conclude
$M_{2,p}(w,r)\leq c\,r^{1-\e_1}$, as stated in Theorem 1. This completes the proof. $\Box$

\medskip
Now let us combine Theorem 1 with a change of variables to treat a general $y\in \RR^n$ and we do not assume $a_{ij}(y)=\delta_{ij}$. We let $B_\e(y)=\{x:|x-y|<\e\}$ and  want to construct a singular solution of
\begin{equation}\label{LZ_y=0}
{\mathcal L}(x,\partial_x)Z_y(x)=0 \quad\hbox{for}\ x\in B_\e(y)\backslash\{y\},
\end{equation}
provided $\e$ is sufficiently small. Since ${\bf A}_y$ is symmetric and positive definite, we can define the symmetric matrix ${\bf A}_y^{-1/2}$. This  enables us to define the function $I_y(r)$ as in \eqref{def:I_y}.
 
\begin{corollary}\label{co:1} 
 For $n\geq 3$, $p\in (1,\infty)$, and fixed $y\in U$, assume that ${\bf A}_y$ is positive definite and the coefficients $a_{ij}$  satisfy
  \eqref{|A_x-A_y|<om}. If $y=0$, we assume \eqref{est:b_k} but otherwise the $b_k$ are bounded on $U$. Then, for $\e$ sufficiently small, there exists a solution $Z_y$ of \eqref{LZ_y=0} in the form
 \begin{subequations}
  \begin{equation} \label{Z_y(x)=}
 Z_y(x)=h_y(|{\bf A}_y^{-1/2}(x-y)|)+v(x),
  \end{equation}
 where
  \begin{equation}
 h_y(r)=\int_r^\e s^{1-n}\,e^{I_y(s)}\,ds\,(1+\zeta_y(r))
  \end{equation}
 with $I_y(r)$ given by \eqref{def:I_y} and $\zeta_y$ satisfying
  \begin{equation}
 M_{2,p}(\zeta_y,r)\leq c\,\max(\om(r),\si(r)),
  \end{equation}
 and $v$ satisfying
  \begin{equation}
 M_{2,p}(v,r;y)\leq c\,r^{2-n}\,e^{I_y(r)}\,\om(r).
  \end{equation}
 \end{subequations}
 Moreover, for any $u\in W^{2,p}_\loc(\overline{B_\e(y)}\backslash\{y\})$ that is a strong solution of ${\mathcal L}(x,\partial_x)u=0$ in 
 $\overline{B_\e(y)}\backslash\{y\}$ and subject to the growth condition
 \[
 M_{2,p}(u,r;y)\leq c\,r^{1-n+\e_0}\quad\hbox{where}\ \e_0>0,
 \]
there exist constants $c_0,c_1$ such that
 \[
 u(x)=c_0\,Z_y(x)+c_1+w(x),
 \]
 where $M_{2,p}(w,r;y)\leq c\, r^{1-\e_1}$ for any $\e_1>0$.
\end{corollary}
\noindent
If we use integration by parts, we can write the solution of Corollary 1 as
\begin{equation}\label{Z_y-asym}
Z_y(x)=\frac{\langle {\bf A}_y^{-1}(x-y),(x-y)\rangle^{\frac{2-n}{2}}}{(n-2)}\,e^{I_y(\sqrt{\langle{\bf A}_y^{-1}(x-y),(x-y)\rangle})}(1+\xi_y(x))
\end{equation}
where $M_{1,\infty}(\xi_y,r;y)\leq c\,\max(\om(r),\si(r))$ for $0<r<\e$. This generalizes \eqref{Z-asym}.

  \medskip\noindent
  {\bf Proof of Corollary 1.} First we consider $y=0$. Since ${\bf A}_0$ is positive definite, we can define a symmetric matrix by ${\bf J}={\bf A}_0^{-1/2}$, so ${\bf J}\, {\bf A}_0 \,{\bf J} ={\bf I}$.  If we introduce a change of variables by $\tilde x={\bf J}\,x$ and new coefficients $\widetilde a_{ij}(\tilde x)$ and $\widetilde b_k(\tilde x)$ by
  \[
  \widetilde {\bf A}_{\tilde x}={\bf J}{\bf A}_x {\bf J} \quad\hbox{and}\quad \widetilde{\bf B}_{\tilde x}={\bf B}_x\,{\bf J},
  \] 
  then $\widetilde {\bf A}_0={\bf I}$ and 
\[
{\mathcal L}(x,\partial_x)=a_{ij}\frac{\partial^2}{\partial x_i\partial x_j}+b_k\frac{\partial}{\partial x_k}
=\widetilde a_{ij}\frac{\partial^2}{\partial \tilde x_i\partial \tilde x_j}+\widetilde b_k\frac{\partial}{\partial \tilde x_k}
=\widetilde{\mathcal L}(\tilde x,\partial_{\tilde x}).
\]
Hence we may apply Theorem 1 in the variables $\tilde x$ to conclude  the existence of a solution $\widetilde Z$ of 
$\widetilde{\mathcal L}(\tilde x,\partial_{\tilde x})\widetilde Z(\tilde x)=0$ for $0<|\tilde x|<\e$ in the form
$\widetilde Z(\tilde x)+\widetilde h(|\tilde x|)+\widetilde v(\tilde x)$ where
\[
\widetilde h(r)=\int_r^\e s^{1-n}\,e^{\widetilde I(s)}\,ds\,(1+\zeta(r))
\]
with $M_{2,p}(\zeta,r)\leq c\,\max(\om(r),\si(r))$, $M_{2,p}(\widetilde v,r)\leq c\, r^{2-n}\,e^{\widetilde I(r)}\,\om(r)$, and
\[
\widetilde I(r)=\frac{1}{|S^{n-1}|}\int_{r<|\tilde x|<\e}\left({\rm tr}(\widetilde{\bf A}_z)-n\,\frac{\langle \widetilde{\bf A}_z z,z\rangle}{|z|^2} 
+\langle \widetilde {\bf B}_z,z\rangle\right)\,\frac{dz}{|z|^n}.
\]
In terms of the original variables, we see that $Z(x)=\widetilde h(|{\bf J}x|)+\widetilde v({\bf J}x)$ satisfies ${\mathcal L}(x,\partial_x)Z=0$ for $0<|{\bf J}x|<\e$, hence for $0<|x|<\e_1$ with $\e_1$ sufficiently small.

Finally, if $y$ is a general point in $U$, then we use the change of variables $\tilde x={\bf J}(x-y)$ with ${\bf J}={\bf A}_y^{-1/2}$
and let $\widetilde{\bf A}_{\tilde x}=(\widetilde{a}_{ij}(\tilde x))={\bf J}{\bf A}_x{\bf J}. $ Since $\tilde x=0$ corresponds to $x=y$, we have $\widetilde a_{ij}(0)=\delta_{ij}$, so we can apply Theorem 1 to $\widetilde{\mathcal L}(\tilde x,\partial_{\tilde x})={\mathcal L}(x,\partial_x)$ to obtain the solution $\widetilde h(|\tilde x|)+\widetilde v(\tilde x)$. We obtain the solution  of \eqref{LZ_y=0} as
\[
Z_y(x)=\widetilde h(|{\bf J}(x-y)|)+\widetilde v({\bf J}(x-y)),
\]
where $\widetilde h(r)$ involves the above $\widetilde I(r)$. To transform $\widetilde I(r)$ to the original variables, replace $\widetilde {\bf A}_{\tilde z}$ by ${\bf A}_z$ and every occurence of $\tilde z$ by $x-y$; we find $\widetilde I$ is of the desired form \eqref{def:I_y}.
Moreover, since $\widetilde v$ satisfies $M_{2,p}(\widetilde v,r)\leq c\,r^{2-n} e^{I(r)}\om(r)$, we find that $v(x)=\widetilde v({\bf J}(x-y))$
satisfies $M_{2,p}(v,r;y)\leq c\,r^{2-n}\,e^{I(r)}\om(r)$, as desired.
$\Box$
 
\section{Finding the constant $C_y$ so that $-{\mathcal L}Z(x)=C_y\delta(x-y)$ }\label{sec:2}

In this section we first prove Theorem 2, then state and prove Corollary 2, which shows how the formulas change when we no longer assume $y=0$ and $a_{ij}(y)=\delta_{ij}$.
As in the proof of Theorem 1, we shall assume \eqref{est:sigma(1)} holds for $\mu$ sufficiently small and work in $B_1$ instead of $B_\e$.
So we want to determine when the singular solution $Z(x)$ of Theorem 1 satisfies
\begin{equation}\label{LZ=delta-B1}
-{\mathcal L}(x,\partial_x)\,Z(x)=C_0\,\delta(x) \quad\hbox{ in $ U=B_1$,}
\end{equation}
 for some constant $C_0$. This was done for ${\mathcal L}_0$ in \cite{MM1}, and many of the arguments here are the same, so we shall be brief; for more details, consult \cite{MM1}. 
 
 From Theorem 1, we obtain estimates on $\partial Z$ and $\partial^2Z$ that show  ${\mathcal L}(x,\partial_x)Z(x)$ can be regularized at $x=0$ as a distribution ${\mathcal F}_0$ on $C^\lambda_0(U)$, the H\"older continuous functions with compact support in $U$. So {\it if} we can define  ${\mathcal L}(x,\partial_x)Z(x)$ as a distribution ${\mathcal F}$, then it must be supported at $x=0$ and only involve $\delta(x)$, not derivatives of $\delta(x)$, i.e.\ satisfy \eqref{LZ=delta-B1} for some constant $C_0$.
 
 The difficulty in defining ${\mathcal L}(x,\partial_x)Z(x)$ as a distribution comes from the lack of regularity of the coefficients, especially $a_{ij}$. In particular, there is no difficulty in defining 2nd-order distributional derivatives of $Z$ by
 \[
 \langle \partial_i\partial_j Z,\phi\rangle = -\int_U \partial_j Z(x)\,\partial_i \phi(x)\,dx \quad\hbox{for}\ \phi\in C^1_0(U),
 \]
since the integral on the right converges. So let us try to define the distribution ${\mathcal L}Z$ by
 \begin{equation}\label{def:LZ}
 \langle {\mathcal L}Z,\phi\rangle
 \!=\!\int_U\! \left( (a_{ij}-\delta_{ij})\partial_i\partial_j Z\phi-\partial_iZ\,\partial_i\phi+b_k\partial_kZ\, \phi\right)dx
  \quad\hbox{for}\ \phi\in C^1_0(U).
 \end{equation}
This is an improper integral due to the singularities in $\partial_i\partial_j Z$ and $b_k\partial_k Z$ at $x=0$, but if the integral converges then we conclude \eqref{LZ=delta-B1} holds and we can compute $C_0$ from
\begin{equation}\label{def:C0}
-C_0=\lim_{\e\to 0} \int_U \left((a_{ij}-\delta_{ij})\,\partial_i\partial_j Z\,\phi_\e-\partial_iZ\,\partial_i\phi_\e+b_k\partial_kZ\,\phi_\e\right)dx,
\end{equation}
where $\phi_\e(|x|)=\chi(|x|/\e)$ with $\chi(r)$ being a smooth cutoff function that is 1 for $0<r<1/4$ and vanishes for $r>1/2$.
(We may assume $\phi(x)=\phi(|x|)$ since we can write $\phi(x)=\phi_0(|x|)+\phi_1(x)$ where $|\phi_1(x)|+|x|\,|\nabla\phi_1(x)|\leq c\,|x|$ for $|x|<1$, which shows that $\langle {\mathcal L}Z,\phi_1\rangle$ is well-defined as an integral and contributes nothing to $C_0$.)

\medskip\noindent
{\bf Proof of Theorem 2.} Recall from the proof of Theorem 1 the decomposition $Z(x)=h(|x|)+v(x)$ in \eqref{Z=h+v}. We first show that $v$ makes no contribution to determining the value of $C_0$. Since $v\in Y$, we have $M_{2,p}(v,r)\leq c\,r^{2-n}\,\om(r)\,e^{I(r)}$.  But $I(r)$ is bounded above in both cases (i) and (ii), so we have
\begin{equation}
M_{2,p}(v,r)\leq c\,r^{2-n}\,\om(r) \quad\hbox{for}\ 0<r<1.
\end{equation}
 Then, as $\e\to 0$, we have\footnote{In the following, integrals involving $\int_{|x|<\e}$ should be interpreted as improper:  $\lim_{\eta\to 0}\int_{\eta<|x|<\e}$.}
 \[
 \left| \int_{|x|<\e} (a_{ij}-\delta_{ij})\,\partial_i\partial_j v\,\phi_\e\,dx\right|\leq c\int_0^\e \om^2(r)\frac{dr}{r}=c\,\sigma(\e)\to 0
 \]
  \[
 \left| \int_{|x|<\e} \partial_i v\, \partial_i\phi_\e\,dx\right|\leq c\,\e^{-1}\int_0^\e \om(r)\,dr\leq c\,\om(\e)\to 0
 \]
 and
 \[
 \left|  \int_{|x|<\e} b_k\,\partial_k v\,\phi_\e\,dx \right|\leq c\int_0^\e \om^2(r)\frac{dr}{r}=c\,\sigma(\e)\to 0.
 \]
 So $v$ makes no contribution to $C_0$.

Now we consider $h(r)$. In fact, from \eqref{h(r)=} we can write $h(r)=h_0(r)+h_1(r)+c$, where $c$ is an arbitrary constant and
 \begin{equation}\label{def:h_0}
 \begin{aligned}
 h_0(r)&=c_*\int_r^1 s^{1-n}\,E_+(s)\,ds, \\
 h_1(r)&=\int_r^1 s^{1-n}\,E_+(s)\int_0^s \rho^{n-1}\,E_-(\rho)\,B[v](\rho)\,d\rho\,ds,
 \end{aligned}
 \end{equation}
 with $E_+$  defined in \eqref{def:E_pm} and $c_*$ chosen as in the proof of Theorem 1 so that $c_*E_+(0)=e^{I(0)}$. (Note that this decomposition of $h(r)$ is slightly different from \eqref{h0,h1}.) Let us show that $h_1$ and $c$ do not contribute to $C_0$.
 We compute
\[
\partial_i h_1=-x_i\,r^{-n}\,E_+(r)\int_0^r \rho^{n-1}\,E_-(\rho)\,B[v](\rho)\,d\rho,
\]
\[
\begin{aligned}
\partial_i\partial_j h_1\!=\!-r^{-n}\,E_+(r)\!\left(\delta_{ij}-(n+R(r))\frac{x_ix_j}{r^2}\right)\!\int_0^r\rho^{n-1}E_-(\rho)B[v](\rho)\,d\rho
\!-\!\frac{x_ix_j}{r^2}B[v](r),
\end{aligned}
\]
and hence
\[
\begin{aligned}
\int_{|x|<\e}\left[ (a_{ij}-\delta_{ij})\partial_i\partial_j h_1  \phi_\e+ b_k\partial_kh_1 \phi_\e-\partial_i h_1\partial_i\phi_\e\right]dx& \\
= \int_0^\e\left[-\frac{R(r)}{r} E_+(r )\chi\left(\frac{r}{\e}\right) +E_+\frac{d}{dr}\chi\left(\frac{r}{\e}\right) \right]&\int_0^r\rho^{n-1}E_-(\rho)B[v](\rho)\,d\rho\,dr\\
-\int_0^\e (\alpha(r)-1)r^{n-1}B[v](r) \chi\left(\frac{r}{\e}\right)\,dr& \\
+ \,\e^{-1}\int_0^\e E_+(r) \int_0^r \rho^{n-1} E_-(\rho)  B[v](\rho)&\,d\rho \, \chi\left(\frac{r}{\e}\right)dr.
\end{aligned}
\]
Let us denote these three terms ${\mathcal I}_1(\e)$, ${\mathcal I}_2(\e)$, ${\mathcal I}_3(\e)$ and estimate them separately.
 First,
\[
\begin{aligned}
{\mathcal I}_1(\e)\!=&\!\int_0^\e \frac{d}{dr}\left[ E_+(r)  \chi\left(\frac{r}{\e}\right)\int_0^r\rho^{n-1} E_-(\rho)B[v](\rho)\,d\rho \right]dr\!-\!\int_0^\e r^{n-1} \chi\left(\frac{r}{\e}\right) B[v](r)\,dr\\
&=E_+(\e)\int_0^e r^{n-1}E_-(r)B[v](r)\,dr-\int_0^\e r^{n-1}B[v](r) \chi\left(\frac{r}{\e}\right)\,dr.
\end{aligned}
\]
But recall from \eqref{est:int-B[v]}  that $\int_0^\e r^{n-1} E_-(r)|B[v](r)|\,dr\leq c\,\si(\e)$ and $I(r)$ is bounded implies $E_\pm(r)$ are bounded above and below by positive constants, so $|{\mathcal I}_1(\e)|\leq c\,\si(\e)\to 0$ as $\e\to 0$. This also implies $|{\mathcal I}_2(\e)|\to 0$ as $\e\to 0$. As for ${\mathcal I}_3(\e)$, we also have
\[
|{\mathcal I}_3(\e)|\leq c\,\e^{-1}\int_0^\e \si(r)\,dr\leq c\,\si(\e)\to 0\quad\hbox{as}\ \e\to 0.
\]
We conclude that $h_1$ does not contribute to $C_0$.
 
 Finally, we consider $h_0$. We can calculate 
 $\partial_i\,h_0=-c_*\,r^{-n}\,E_+(r)\,x_i$,  and 
\[
\del_i\del_j h_0=-c_*\,r^{-n}\,E_+(r)\left(\de_{ij}-n\frac{x_ix_j}{r^2}-\frac{x_ix_j}{r^2}R(r)\right).
\]
It is easy to verify that
\[
(a_{ij}-\de_{ij})\del_i\del_j h_0 = -c_*\, r^{-n}\, E_+(r) \left(a_{ii}-n\frac{a_{ij}x_ix_j}{r^2}-\left(\frac{a_{ij}x_ix_j}{r^2}-1\right)R(r)\right).
\]
Notice that
\[
\int_{|x|<\e} r^{-n}\, E_+(r) \left(a_{ii}-n\frac{a_{ij}x_ix_j}{r^2}-\left(\frac{a_{ij}x_ix_j}{r^2}-1\right)R(r)\right)\phi_\e(|x|)\,dx
\]
\[
=|S^{n-1}|\int_0^\e E_+(r) \frac{\alpha_n(r)-n\alpha(r)-(\alpha(r)-1)R(r)}{r}\chi\left(\frac{r}{\e}\right)\,dr
\]
\[
=|S^{n-1}|\int_0^\e E_+(r)\frac{R(r)-\beta(r)}{r}\chi\left(\frac{r}{\e}\right)\,dr,
\]
since $\alpha(r)\,R(r)=\alpha_n(r)+\beta(r)-n\alpha(r)$. Similarly, we can verify
\[
\int_{|x|<\e} b_k\partial_k h_0\,\phi_\e\,dx=-c_*\,|S^{n-1}|\int_0^\e E_+(r)\,\frac{\beta(r)}{r}\,\chi\left(\frac{r}{\e}\right)\,dr,
\]
so
 \[
  \langle {\mathcal L}h_0,\phi_\e\rangle=c_*\,|S^{n-1}|\int_0^\e \left(-E_+(r)\frac{R(r)}{r}\chi\left(\frac{r}{\e}\right)+E_+(r)\frac{d}{dr}\left[\chi\left(\frac{r}{\e}\right)\right]\right)\,dr.
 \]
Hence
 \begin{equation}
 -\langle {\mathcal L}h_0,\phi_\e\rangle=-c_*|S^{n-1}|\,\int_0^\e  \frac{d}{dr} \left[ E_+(r) \chi\left(\frac{r}{\e}\right)\right]\,dr=c_*\,|S^{n-1}|\,E_+(0).
 \end{equation}
 Letting $\e\to 0$, we conclude $Z$ satisfies \eqref{LZ=delta-B1} with $C_0=|S^{n-1}| e^{I(0)}$.
 $\Box$
 
 \medskip
 Now let us consider a general fixed $y\in U$ and try to solve
 \begin{equation}\label{LZ_y=C0*delta}
 -{\mathcal L}(x,\partial_x)Z_y(x)=C_y\,\delta(x-y) \quad\hbox{for}\ x\in B_1(y)
 \end{equation}
 for some constant $C_y$. We replace \eqref{def:LZ} with
  \begin{equation}\label{def:LZ_y}
 \langle {\mathcal L}Z_y,\phi\rangle
 \!=\!\int_U\! \left( (a_{ij}-a_{ij}(y))\partial_i\partial_j Z_y\phi-a_{ij}(y)\,\partial_iZ_y\,\partial_j\phi+b_k\partial_kZ_y\, \phi\right)dx.
 \end{equation}
\begin{corollary}\label{co:2} 
Suppose the conditions of Corollary 1 hold and $Z_y$ is the singular solution found there with $I_y(r)$ given by \eqref{I_y=}. 
         \begin{itemize}
\item[(i)] If $I_y(0)=\lim_{r\to 0}I_y(r)$ exists and is finite, then
we can solve \eqref{LZ_y=C0*delta} in $B_1(y)$ with $C_y=|S^{n-1}|\,\sqrt{{\rm det}{\bf A}_y}\,e^{I_y(0)}$.
\item[(ii)] If $I_y(r)\to -\infty$ as $r\to 0$, then solving  \eqref{LZ=C0*delta} in $B_1(y)$ yields $C_0=0$, and so $Z$ solves $-{\mathcal L}(x,\partial_x) Z_y(x)=0$ in $B_1(y)$, despite its singularity at $x=y$.
\end{itemize}
\end{corollary}
\noindent
If $I_y(0)$ exists and is finite, we see from \eqref{Z_y-asym} that the solution $F_y$ of $-{\mathcal L}(x,\partial_x)F_y(x)=\delta(x-y)$ in $B_\e(y)$ has the asymptotic behavior
\begin{equation}\label{F_y-asym}
F_y(x)=\frac{\langle {\bf A}_y^{-1}(x-y),(x-y)\rangle^{\frac{2-n}{2}}}{(n-2)|S^{n-1}|\sqrt{{\rm det}{\bf A}_y}}\,(1+o(1))
\quad\hbox{as}\ x\to y.
\end{equation}
This generalizes \eqref{F-asym} and establishes \eqref{F(x,y)=}.

 \medskip\noindent
  {\bf Proof of Corollary 2.} We only need to show that
  \[
  \langle -{\mathcal L}Z_y,\phi\rangle = |S^{n-1}|\,({\rm det}{\bf A}_y)^{-1/2}\,e^{I_y(0)}\phi(y)
  \]
  for some $\phi\in C^\infty_0(B_{\e_y}(y))$. Let us recall the change of coordinates used in the proof of Theorem 1, namely $\tilde x={\bf J}(x-y)$ where ${\bf J}={\bf A}_y^{-1/2}$, and let $\widetilde\phi(\tilde x)=\phi(x)$. Then
  \[
  -\int {\mathcal L}(x,\partial_x) Z_y(x)\phi(x)\,dx=-({\det}{\bf A}_y)^{1/2}\int \widetilde{\mathcal L}(\tilde x,\partial_{\tilde x})\widetilde Z_0(\tilde x)\widetilde\phi(\tilde x)\,d\tilde x.
  \]
  But Theorem 2 implies $-\langle \widetilde{\mathcal L}\widetilde Z_0,\widetilde\phi\rangle=|S^{n-1}|\,e^{I_y(0)}\widetilde\phi(0)$.
  Since $\widetilde\phi(0)=\phi(y)$, we obtain the desired result. $\Box$
  
\section{Constructing the fundamental solution}\label{sec:3}

\noindent
  {\bf Proof of Theorem 3.} For each $y\in U$, denote the $\e$ in Corollary 1 by $\e_y$, and  use Corollary 2 to calculate $C_y$, which is positive since $I_y(0)$ is finite. We conclude that
  $-{\mathcal L}(x,\partial_x) Z_y(x)/C_y=\delta(x-y)$ for all $x,y\in U$ with $|x-y|<\e$.  For fixed $y\in U$, let $\eta_y(r)$ be a smooth cutoff function satisfying $\eta_y(r)=1$ for sufficiently small $r$ and define
  \begin{equation}\label{def:F}
  F(x,y)=\eta_y(|x-y|) Z_y(x)/C_y +v(x,y),
  \end{equation}
 where $v(x,y)$ is to be determined. But if we apply $-{\mathcal L}(x,\partial_x)$ we obtain
 \[
 -{\mathcal L}(x,\partial_x)F(x,y)=\delta(x-y)+\psi(x,y)-{\mathcal L}(x,\partial_x)v(x,y),
 \]
  where $\psi(\cdot,y)\in L^p(U)$ for $1<p<\infty$ and vanishes near $y$. We may assume that $U$ has a smooth boundary (else we can embed it in such a bounded domain and extend the coefficients, see \cite{Mi}),
  so we may find $v(x,y)$ by solving the Dirichlet problem 
  (for fixed $y\in U$):
  \begin{equation}\label{DirichletProb}
  \begin{aligned}
   -{\mathcal L}(x,\partial_x)v(x,y)&=\psi(x,y) \quad\hbox{for}\ x\in U,\\
   v(x,y)&=0 \quad\hbox{for}\ x\in\partial U.
   \end{aligned}
  \end{equation}
  When $b_k\in L^\infty(U)$,
 it is well-known (cf.\ Theorem 9.15 in \cite{GT}) that \eqref{DirichletProb} has a unique solution $v(\cdot,y)\in W^{2,p}(U)\cap W^{1,p}_0(U)$
 for $1<p<\infty$; associated with this unique solvability is the \`a priori inequality $\|v\|_{W^{2,p}(U)}\leq c\,\|{\mathcal L}v\|_{L^p(U)}$ for $v\in W^{2,p}(U)\cap W_0^{1,p}(U)$. But our assumption \eqref{est:b_k} enables us to write ${\mathcal L}={\mathcal L}_1+{\mathcal L}_2$, where
 ${\mathcal L}_1$ has bounded coefficients in $U$ and  ${\mathcal L}_2=\tilde b_k\partial_k$ has coefficients supported in a 
 ball $B_\sigma$. For any $\e>0$, we can take $\sigma$ sufficiently small that  $\sup_{x\in B_\sigma} |\tilde b_k(x)|\leq \e\,|x|^{-1}$.
 For $1<p<n$, by Hardy's inequality we have $\|{\mathcal L}_2w\|_{L^p}\leq \e\,\|w\|_{W^{2,p}}$ for $w\in W^{2,p}$ with support in $B_\sigma$. If we take $\e$ sufficiently small, we can arrange that the \`a priori inequality $\|v\|_{W^{2,p}(U)}\leq c\,\|{\mathcal L}v\|_{L^p(U)}$ holds for $v\in W^{2,p}(U)\cap W_0^{1,p}(U)$. We conclude that \eqref{DirichletProb} admits a unique solution $v(\cdot,y)\in W^{2,p}(U)\cap W_0^{1,p}(U)$ for $1<p<n$. Taking $p\in (n/2,n)$, we see that $v(x,y)$ is continuous in $x\in U$.
  
  To confirm that $F(x,y)$ satisfies \eqref{F(x,y)=}, we know from \eqref{def:F} and \eqref{Z_y-asym} that
  \[
  F(x,y)=\frac{\langle {\bf A}_y^{-1}(x-y),(x-y)\rangle^{\frac{2-n}{2}}}{(n-2)|S^{n-1}|\sqrt{\det {\bf A}_y}}\,e^{I_y(\sqrt{\langle{\bf A}_y^{-1}(x-y),(x-y)\rangle})-I_y(0)}(1+\xi_y(x))
  \]
where $M_{1,\infty}(\xi_y,r;y)\leq c\,\max(\om(r),\si(r))$ for $0<r<\e$. Since the exponential term tends to 1 as $|x-y|\to 0$, we obtain 
 \eqref{F(x,y)=}. $\Box$

\begin{appendix}
\section{$L^p$-mean estimates for convolutions}\label{App1}

In \cite{MM1}, we proved $L^p$-mean estimates for distribution solutions of 
\begin{equation}\label{Delta u=f}
\Lap u =f \quad \hbox{in}\ \RR^n\backslash\{0\},
\end{equation}
when $f$ has zero spherical mean, i.e.\ $\overline{f}(r)=0$. Let $K$ denote convolution by the fundamental solution $\Gamma$ of the Laplacian.
The following appears as Corollary 1 in Section 1 of \cite{MM1}.
\begin{proposition}\label{pr:convolution}
Suppose $n\geq 2$, $p\in (1,\infty)$, and $f\in L^p_{\loc}(\RR^n\backslash\{0\})$ satisfies
\begin{equation}
\overline{f}(r)=0 \quad\hbox{and}\quad \int_{|x|<1}|x|\,|f(x)|\,dx+\int_{|x|>1} |x|^{1-n}|f(x)|\,dx <\infty.
\end{equation}
Then $u=Kf=\Gamma\star f$ is a distribution solution of \eqref{Delta u=f} that satisfies
\begin{equation}
M_{2,p}(Kf,r)\leq c\left( r^{1-n}\int_0^r M_p(r,\rho)\,\rho^n\,d\rho +r\int_r^\infty M_p(f,\rho)\,d\rho\right).
\end{equation}
\end{proposition}

\section{The Laplacian on weighted Sobolev spaces}\label{App2}

Many authors have studied the mapping properties of elliptic operators such as the Laplacian on weighted Sobolev spaces on $\RR^n$ and other noncompact manifolds with conical or cylindrical ends: cf.\ \cite{LM}, \cite{MP}, \cite{Mc0}. We will recall some of these results for punctured Euclidean space $\RR^n_o:=\RR^n\backslash\{0\}$. Since we are mostly concerned in this paper with singularities at $x=0$, let us first investigate the weighted $L^p$-norm on the punctured unit ball $B_o:=B_1\backslash\{0\}$. For $\delta\in\RR$ and $1<p<\infty$, define the Banach space $L^p_\delta(B_o)$ by the norm
\begin{equation}\label{def:L^p_delta}
\|u\|^p_{L^p_{\delta}(B_o)}:=\int_{|x|<1} |x|^{\delta p}\,|u(x)|^p\,dx.
\end{equation}
For example, the constants are in $L^p_{\delta}(B_o)$ if and only if $\delta>-n/p$.
To compare \eqref{def:L^p_delta} with $L^p$-means, it is easy to see (cf.\ \cite{MM1}) that $M(u,r)\leq c\,r^\e$ for $0<r<1$ implies that $u\in L_\delta^p(B_o)$ provided $\e+\delta>-n/p$, and conversely $u\in L_\delta^p(B_o)$ implies $M_p(u,r)\leq c_\e\,r^\e$ for $0<r<1$ if $\e=-\delta-n/p$.

Now we introduce the weighted $L^p$-norm for functions on $\RR^n_o$ with separate weights at the origin and at infinity. For $\delta_0,\delta_1\in\RR$, define
\begin{equation}
\begin{aligned}
\|u\|^p_{L^p_{\delta_0,\delta_1}}:=\int_{|x|<1} |x|^{\delta_0 p}|u(x)|^p\,dx+\int_{|x|>1} |x|^{\delta_1 p}|u(x)|^p\,dx.
\end{aligned}
\end{equation}
We then define the weighted Sobolev space $W^{2,p}_{\delta_o,\delta_1}(\RR^n_o)$ to be functions $u\in W^{2,p}_{\loc}(\RR^n_o)$ for which
\begin{equation}
\|u\|_{W^{2,p}_{\delta_o,\delta_1}}=\sum_{|\alpha|\leq 2} \| |x|^{|\alpha|}\partial^\alpha u\|_{L^p_{\delta_1,\delta_2}}
\end{equation}
is finite. It is clear that 
\begin{equation}\label{Lap:W->L}
\Lap: W^{2,p}_{\delta_0,\delta_1}(\RR^n_o)\to L^{p}_{\delta_0+2,\delta_1+2}(\RR^n_o)
\end{equation}
is a bounded linear operator, and (using the analysis of  \cite{LM}, \cite{MP}, \cite{Mc0}, for example) it can be shown that
\eqref{Lap:W->L} is an isomorphism for $-n/p<\delta_0,\delta_1<-2+n/p'$, where $p'=p/(p-1)$. (Since $n\geq 3$, such $\delta_0,\delta_1$ exist.)
Moreover, provided $\delta_0,\delta_1$ do not take the values $-n/p-k$ or $-2+n/p'+k$ where $k$ is a nonnegative integer, then \eqref{Lap:W->L} is a Fredholm operator whose nullspace and cokernel are easily described in terms of harmonic polynomials.
Since we are principally interested in the behavior of functions at the origin, we will fix 
\begin{equation}\label{delta1-range}
-n/p< \delta_1<-2+n/p'
\end{equation}
and allow $\delta_0$ to vary. We only require a small range of  values for $\delta_0$.

\begin{proposition}
Assume $1<p<\infty$ and \eqref{delta1-range}. Then the map \eqref{Lap:W->L} is
\begin{itemize}
\item[(a)] an isomorphism for $-n/p<\delta_0<-2+n/p'$;
\item[(b)] surjective with nullspace spanned by $|x|^{2-n}$ if  $-2+n/p'<\delta_0<-1+n/p'$;
\item[(c)] injective with cokernel spanned by $1$ if $-1-n/p<\delta_0<-n/p$.
\end{itemize}
\end{proposition}

\end{appendix}

 \bigskip\noindent
{\bf Acknowledgement}: This paper has been supported by the RUDN University Strategic
Academic Leadership Program.


\end{document}